\documentclass[10pt]{article}

\usepackage{geometry,setspace}                % See geometry.pdf to learn the layout options. There are lots.
\usepackage{graphicx}
\usepackage{amssymb}
\usepackage{epstopdf}
\usepackage{latexsym, amssymb, amsmath, amsthm, graphics}
\usepackage{url}
\usepackage{eucal}

\title{ \bf The tame automorphism group \\ in two variables \\ over basic Artinian rings} 
\author{Joost Berson}

    \date{}
    
\begin{document}
    \maketitle

 \input amssym.def
    \theoremstyle{definition}
    \newtheorem{rema}{Remark}[section]
    \newtheorem{questions}[rema]{Questions}
    \newtheorem{assertion}[rema]{Assertion}
    \theoremstyle{plain}
    \newtheorem{propo}[rema]{Proposition}
    \newtheorem{theo}[rema]{Theorem}
    \newtheorem{conj}[rema]{Conjecture}
    \newtheorem{quest}[rema]{Question}
    \theoremstyle{definition}
    \newtheorem{defi}[rema]{Definition}
    \theoremstyle{plain}
    \newtheorem{lemma}[rema]{Lemma}
    \newtheorem{corol}[rema]{Corollary}
    \newtheorem{exam}[rema]{Example}
    \newtheorem{rmk}[rema]{Remark}
    \newcommand{\del}{\triangledown}
    \newcommand{\nno}{\nonumber}
    \newcommand{\lbar}{\big\vert}
    \newcommand{\mbar}{\mbox{\large $\vert$}}
    \newcommand{\p}{\partial}
    \newcommand{\dps}{\displaystyle}
    \newcommand{\bra}{\langle}
    \newcommand{\ket}{\rangle}
    \newcommand{\kr}{\mbox{\rm Ker}\ }
    \newcommand{\res}{\mbox{\rm Res}}
    \renewcommand{\hom}{\mbox{\rm Hom}}
    \newcommand{\pf}{{\it Proof:}\hspace{2ex}}
    \newcommand{\epf}{\hspace{2em}$\Box$}
    \newcommand{\epfv}{\hspace{1em}$\Box$\vspace{1em}}
    \newcommand{\nord}{\mbox{\scriptsize ${\circ\atop\circ}$}}
    \newcommand{\wt}{\mbox{\rm wt}\ }
    \newcommand{\clr}{\mbox{\rm clr}\ }
    \newcommand{\ideg}{\mbox{\rm Ideg}\ }
    \newcommand{\gC}{{\mathfrak g}_{\mathbb C}}
    \newcommand{\hatC}{\widehat {\mathbb C}}
    \newcommand{\Z}{{\mathbb Z}}
    \newcommand{\bQ}{{\mathbb Q}}
    \newcommand{\bR}{{\mathbb R}}
    \newcommand{\N}{{\mathbb N}}
    \newcommand{\bT}{{\mathbb T}}
    \newcommand{\fg}{{\mathfrak g}}
    \newcommand{\fgC}{{\mathfrak g}_{\bC}}
    \newcommand{\cD}{\mathcal D}
    \newcommand{\cP}{\mathcal P}
    \newcommand{\cC}{\mathcal C}
    \newcommand{\cS}{\mathcal S}
    \newcommand{\EGC}{{\cal E}(\GC)}
    \newcommand{\cLGC}{\widetilde{L}_{an}\GC}
    \newcommand{\LGC}{{L}_{an}\GC}
    \newcommand{\BQ}{\begin{eqnarray}}
    \newcommand{\EQ}{\end{eqnarray}}
    \newcommand{\BQn}{\begin{eqnarray*}}
    \newcommand{\EQn}{\end{eqnarray*}}
    \newcommand{\wtilde}{\widetilde}
    \newcommand{\Hol}{\mbox{Hol}}
    \newcommand{\Hom}{\mbox{Hom}}
    \newcommand{\poly}{polynomial }
    \newcommand{\polys}{polynomials }
    \newcommand{\pz}{\frac{\p}{\p z}}
    \newcommand{\pzi}{\frac{\p}{\p z_i}}
    \newcommand{\edge}{\text{\raisebox{2.75pt}{\makebox[20pt][s]{\hrulefill}}}}
    \newcommand{\halfedge}{\text{\raisebox{2.75pt}{\makebox[10pt][s]{\hrulefill}}}}
    \newcommand{\n}{\notag}
    \newcommand{\C}{\mathbb C}
    \newcommand{\A}{\mathcal A}
    \newcommand{\Q}{\mathbb Q}
    \newcommand{\X}{X_1,\ldots,X_n}
    \newcommand{\pa}{\partial}
    \newcommand{\D}{\text{D}}
    \newcommand{\Del}{\text{\raisebox{2.1pt}{$\bigtriangledown$}}}
    \newcommand{\La}{\triangle}
    \newcommand{\Hess}{\text{\rm Hess}}
\newcommand{\GA}{\text{GA}}
\newcommand{\SA}{\text{SA}}
\newcommand{\MA}{\text{MA}}
\newcommand{\GL}{\text{GL}}
\newcommand{\SL}{\text{SL}}
\newcommand{\GE}{\text{GE}}
\newcommand{\Tr}{\text{Tr}}
\newcommand{\Af}{\text{Af}}
\newcommand{\Bf}{\text{Bf}}
\newcommand{\W}{\text{W}}
\newcommand{\EA}{\text{EA}}
\newcommand{\BA}{\text{BA}}
\newcommand{\TA}{\text{TA}}
\newcommand{\E}{\text{E}}
\newcommand{\B}{\text{B}}
\newcommand{\T}{\text{T}}
\newcommand{\Di}{\text{D}}
\newcommand{\J}{\text{J}}
\newcommand{\Ker}{\text{Ker}}
\newcommand{\degr}{\text{deg}\,}
\newcommand{\supp}{\text{supp}\,}
\newcommand{\F}{\mathbb F}
\newcommand{\ideaala}{\mathfrak{a}}
\newcommand{\ideaalb}{\mathfrak{b}}
\newcommand{\ideaalc}{\mathfrak{c}}
\newcommand{\m}{\mathfrak{m}}

\renewcommand{\thefootnote}{}

\abstract{
\noindent In a recent paper it has been established that over an 
Artinian ring $R$ all two-dimensional polynomial automorphisms having Jacobian 
determinant one are tame if $R$ is a $\Q$-algebra.  This is a generalization 
of the famous Jung-Van der Kulk Theorem, which deals with the case that $R$ 
is a field (of any characteristic). Here we will show that for tameness over 
an Artinian ring, the $\Q$-algebra assumption is really needed: we will give, 
for local Artinian rings with square-zero principal maximal ideal, 
a complete description of the tame automorphism subgroup. This will lead to 
an example of a non-tame automorphism, for any characteristic $p>0$.\\
}


\footnote{Funded by a Free Competition grant from the Netherlands Organisation for Scientific Research (NWO)\\}
\footnote{Joost Berson, Radboud University, Faculty of Science, P.O. Box 9010, 6500 GL Nijmegen, The Netherlands,}
\footnote{j.berson@science.ru.nl}

\section{Introduction} All two-dimensional polynomial automorphisms over a field are tame,
as stated in the famous theorem by Jung and Van der Kulk (\cite{Jung},\cite{vdK}). 
For fields of characteristic zero this was proved by Jung, and Van der Kulk generalized it 
to arbitrary characteristic.  As is well-known, the statement fails to be true over a domain $R$ 
which is not a field. The most common example of a non-tame automorphism is the one 
by Nagata (\cite{Nag}), which is defined over $R=k[Z]$ (a univariate polynomial ring), 
but can be transformed into an example over any domain which is not a field. 
For any domain $R$, \cite[Corollary~5.1.6]{Essen} even yields an algorithm to decide 
whether or not an automorphism in two variables over $R$ is tame.
To continue with the description of the tame automorphism groups over commutative rings 
in general, it is very convenient to start with Artinian rings. Namely, when it is clear 
which automorphisms in two variables over Artinian rings are tame, we can use this information 
to decribe the automorphisms over rings with higher Krull dimension, by lifting 
the former automorphisms (see for example Theorem~\ref{modnil}). Moreover, the problem of 
describing the structure of the automorphism group over a general Artinian ring can be reduced 
to the case of a local Artinian ring (as will be explained in section~\ref{specialcase}).

One of the main results of the recent paper \cite{BEW08} by Van den Essen, Wright and the author 
is the fact that, over an Artinian ring $R$, all two-dimensional automorphisms with 
Jacobian determinant one are tame in case $R$ is a $\Q$-algebra. This is a generalization of 
Jung's Theorem. We show that in a non-$\Q$-algebra setting, tameness is not guaranteed. 
In fact, for every characteristic $p>0$ we give an example of a non-tame automorphism over 
a local Artinian ring having that characteristic. To show that these automorphisms are not tame, 
we provide a description of the structure of the tame automorphism groups over local Artinian rings 
of the most basic type to be found: the ones with square-zero principal maximal ideal.

This paper is set up as follows: In the next section we introduce the general automorphism 
group and its best-known subgroups. We describe classic results, explain which questions 
are still unanswered and in what way this paper contributes to the development of the theory 
on polynomial automorphism groups. In section~3, we review one of the results of the recent 
paper \cite{BEW08}, saying that over an Artinian $\Q$-algebra $R$, any two-dimensional 
automorphism is tame, provided that the Jacobian determinant is equal to one 
(Theorem~\ref{dimzerotame}). The preparations for this result, which will be done in that section, 
are also important for the remainder of this paper: most techniques also work in the non-$\Q$-algebra setting. 
Lemma~\ref{linearcomb} is in fact the only tool that requires a $\Q$-algebra. 
Section~4 examines the structure of the elementary automorphism subgroup $\EA_2(R)$ 
for rings $R$ of the form $R=A[T]/(T^2)$, where $A$ is another ring. It essentially 
reduces the description of the elementary subgroup over $R$ to the description of 
the elementary subgroup over $A$. This result can immediately be applied to the case 
of local Artinian rings with square-zero principal maximal ideal. This is done in the 
last section. It yields an example for every prime number $p$ of a non-tame automorphism 
in two variables over $\F_p[T]/(T^2)$.

\section{Automorphism subgroups and their relations}

In this paper, every ring is assumed to be commutative and to have an identity element.
We will restrict ourselves to polynomial rings in two variables over a ring $R$, denoted as $R[X,Y]$. 
This section describes the usual subgroups of the general polynomial automorphism group, 
and what is already known about how they are related.

A \textit{polynomial map over $R$} is an ordered pair $(F,G)$ of polynomials of $R[X,Y]$.
We can view polynomial maps as maps $R^2\rightarrow R^2$, defined by $(x,y)\mapsto(F(x,y),G(x,y))$, 
but also as $R$-endomorphisms $R[X,Y]\rightarrow R[X,Y]$, given by the substitution $h(X,Y)\mapsto h(F,G)$. 
$F$ and $G$ are called the \textit{coordinates} of $(F,G)$.

In the usual notation, the composition of two polynomial maps $(F_1,G_1)$ and $(F_2,G_2)$ 
is defined as $(F_1,G_1)\circ(F_2,G_2)=(F_1(F_2,G_2),G_1(F_2,G_2))$. The map $(F_1,G_1)$ 
is called an \textit{invertible polynomial map} or an \textit{automorphism} if there exist 
a polynomial map $(F_2,G_2)$ with $(F_1,G_1)\circ(F_2,G_2)=(F_2,G_2)\circ(F_1,G_1)=(X,Y)$ (the identity map). 
The automorphisms form a group, $\GA_2(R)$.

We write $\J\varphi$ for the Jacobian matrix of an automorphism $\varphi$. 
By the chain rule, for any automorphism $\varphi$ we have $\J\varphi \in \GL_2(R)$, 
whence $|\J\varphi| \in R[X]^*$. 
(Throughout this paper, the operator $|\,\,\,|$ takes the determinant of a matrix.)

Here is an overview of the usual subgroups of $\GA_2(R)$:
\begin{enumerate}
\item $\SA_2(R)$, the {\it special automorphism group}, is the subgroup of all $\varphi$ for which $|J\varphi|=1$. 
\item The group $\GL_2(R)$ of invertible matrices is usually viewed as a subgroup of $\GA_2(R)$.
\item $\EA_2(R)$ is the subgroup generated by the elementary automorphisms.  An {\it elementary} automorphism is one of the form 
$(X+f(Y),Y)$ or $(X,Y+f(X))$ for some univariate polynomial $f$.\\
Note that $\EA_2(R)\subseteq \SA_2(R)$.
\item $\TA_2(R)$, the group of {\it tame} automorphisms, is the subgroup generated by 
$\GL_2(R)$ and $\EA_2(R)$. 
\end{enumerate}
In the case of a field we have the following classic theorem, 
which was already mentioned in the introduction.

\begin{theo}[Jung \cite{Jung} - Van der Kulk \cite{vdK}] \label{JvdK}
For any field $k$ we have $\TA_2(k)=\GA_2(k)$.
\end{theo}

\noindent So over a field the only examples of polynomial automorphisms are the tame ones. 
However, this doesn't hold for a domain which is not field. But there exists an 
algorithm to decide whether or not an automorphism in two variables over a domain $R$ is tame 
in \cite[Corollary~5.1.6]{Essen}). This algorithm can be used to show that any non-unit 
$r \in R\backslash\{0\}$ produces a non-tame automorphism, namely 
$$\left(X-2Y(rX+Y^2)-r(rX+Y^2)^2,Y+r(rX+Y^2)\right)$$
For a polynomial ring $R=k[Z]$ and $r=Z$, $k$ a field, this is Nagata's famous example (\cite{Nag}). 
But for a general commutative ring little is known about which automorphisms in $\GA_2(R)$ are tame. 
This paper is meant to extend our knowledge on this subject.

If $R\to S$ is a surjective ring homomorphism, then the induced group homomorphism $\EA_2(R)\to\EA_2(S)$ is also surjective. 
Note that this fails to hold for $\TA_2(R)\to\TA_2(S)$, because of the following: if $M \in \GL_2(S)$, then 
there doesn't necessarily exist an $N \in \GL_2(R)$ such that $N\mapsto M$. 
This is why tame automorphisms appearing in this paper are usually elements of $\EA_2(-)$: 
we can lift these to automorphisms over rings with higher Krull dimension. 
Therefore, we would like to know the connection between $\TA_2(R)$ and $\EA_2(R)$. 
The following lemma (a special version of \cite[Proposition~3.20]{BEW08}) and corollary 
decribe this connection, which applies to most coefficient rings considered in this paper. 

In the following, $\SL_2(R)$ denotes the group of all matrices with determinant one, $\D_2(R)$ is the group 
of all invertible diagonal matrices, and $\E_2(R)$ is the group generated by all elementary matrices.

\begin{lemma} \label{primtam}
If $R$ is a ring for which $\SL_2(R)=\E_2(R)$, then $\TA_2(R)\cap\SA_2(R)=\EA_2(R)$. 
The hypothesis holds when $R$ is a local ring.
\end{lemma}

\begin{proof}
From $\GL_2(R)=\langle\SL_2(R),\D_2(R)\rangle=\langle\E_2(R),\D_2(R)\rangle\subseteq\langle\EA_2(R),\D_2(R)\rangle$ 
we get that $\TA_2(R)=\langle\GL_2(R),\EA_2(R)\rangle\subseteq\langle\D_2(R),\EA_2(R)\rangle$, 
whence $\TA_2(R)=\langle\D_2(R),\EA_2(R)\rangle$. Since one can then 
readily verify that $\EA_2(R) \triangleleft \TA_2(R)$, this implies that $\TA_2(R)=\D_2(R)\EA_2(R)$. 
But then $\TA_2(R)\cap\SA_2(R) = (\D_2(R)\cap\SA_2(R))\EA_2(R) = (\D_2(R)\cap\SL_2(R))\EA_2(R) = (\D_2(R)\cap\E_2(R))\EA_2(R) 
\subseteq \EA_2(R)$. This proves the first statement. 

For the second statement, consider an element $M$ of $\SL_2(R)$, where $R$ is local. 
Since $\det(M)=1$, there must at least be one entry of $M$ which is in $R^*$. 
We can use this entry to clear the other entries of the row and column to which 
this entry belongs, through multiplication by 2 elementary matrices. 
If the resulting matrix isn't diagonal, we can make it so by multiplying it with 
the matrix 
$$
\left(\begin{array}{cc} 0 & 1 \\ -1 & 0
\end{array}\right)=\left(\begin{array}{cc} 1 & 1 \\ 0 & 1
\end{array}\right) \left(\begin{array}{cc} 1 & 0 \\ -1 & 1
\end{array}\right) \left(\begin{array}{cc} 1 & 1 \\ 0 & 1
\end{array}\right) \in \E_2(R)
$$
Hence, we may assume that the resulting matrix is diagonal, and since it is 
still an element of $\SL_2(R)$, we can use the fact that, for any ring $R$ and any $a \in R^*\!$, 
$$
\left(\begin{array}{cc} a & 0 \\ 0 & a^{-1} \end{array}\right)=
\left(\begin{array}{cc} 1 & a \\ 0 & 1 \end{array}\right)
\left(\begin{array}{cc} 1 & 0 \\ -a^{-1} & 1 \end{array}\right)
\left(\begin{array}{cc} 1 & a-1 \\ 0 & 1 \end{array}\right)
\left(\begin{array}{cc} 1 & 0 \\ 1 & 1 \end{array}\right)
\left(\begin{array}{cc} 1 & -1 \\ 0 & 1 \end{array}\right)
 \in \E_2(R)
$$
\end{proof}

\begin{rema}
If $R$ is not assumed to be local, then the hypothesis of Lemma~\ref{primtam} 
still holds if $R$ has a special structure, \emph{e.g.} when $R$ is a Euclidean domain 
(the proof of this well-known fact is very much like the proof of the second part of Lemma~\ref{primtam}). 
It is important to note, however, that not all Principal Ideal Domains have this 
property. Let $R$ be the ring of integers of $\Q(\sqrt{-19})$, one of the finitely many imaginary quadratic number fields 
of which the ring of integers is a Principal Ideal Domain, by the Stark-Heegner Theorem (\cite{stark}, \cite{heegner}).
In \cite[Theorem~6.1]{Cohn} it was shown that, if $\alpha:=\genfrac{}{}{}{1}{1}{2}+\genfrac{}{}{}{1}{1}{2}\sqrt{-19}$, 
then the following matrix is in $\SL_2(R)$, but not in $\E_2(R)$:
$$
\left(\begin{array}{cc} 3-\alpha & 2+\alpha \\ -3-2\alpha & 5-2\alpha \end{array}\right) 
$$
\end{rema}

\vspace{0.5cm}

\begin{corol} \label{TAandSA}
For any ring $R$ we have the following: if $\EA_2(R)=\SA_2(R)$, then $\TA_2(R)=\{\varphi \in \GA_2(R)\,:\,|\J\varphi| \in R^*\}$. The reverse holds 
if $R$ is any ring for which $\SL_2(R)=\E_2(R)$.
\end{corol}

\begin{proof}
For the first statement, let $\varphi \in \GA_2(R)$ with $|\J\varphi| \in R^*$ 
(since $R$ is reduced). Then there exists an $\alpha \in \GL_2(R)$ such that $\alpha\varphi \in \SA_2(R)=\EA_2(R)$. 
Thus, $\varphi \in \TA_2(R)$.

The second statement follows directly from Lemma~\ref{primtam}.
\end{proof}

\section{The Artinian $\Q$-algebra result}

\noindent Throughout this section (except for Lemma~\ref{sumcomp}), we assume that $R$ is a $\Q$-algebra. We will 
restate and give a quick proof of one of the results from \cite{BEW08}: for an Artinian 
$\Q$-algebra $R$, every special automorphism in two variables over $R$ is tame 
(Theorem~\ref{dimzerotame}). The fact that this is also true for any reduced Artinian ring 
($\Q$-algebra or not) had already been observed in \cite[Corollary~0.6]{Nag} and \cite[Proposition~3.10]{stabtam2}. 
One of the basic tools of Theorem~\ref{dimzerotame} is Lemma~\ref{sumcomp}. 
This lemma is also useful for the general (non-$\Q$-algebra) case in the subsequent sections.
Lemma~\ref{linearcomb} is taken from \cite{BEW08}, and its statement also appeared 
in \cite[\S~5.2, Exercise~7]{Essen}. It is the only ingredient of Theorem~\ref{dimzerotame} 
that requires $R$ to be a $\Q$-algebra.

\vspace{0.1cm}

\begin{lemma} \label{linearcomb}
Every monomial $X^nY^m$ in $R[X,Y]$ can be written as
a $\Q$-linear combination of polynomials of the form $(X+aY)^{n+m}$,
with $a \in \Q$.
\end{lemma}

\vspace{0.1cm}

\noindent The following lemma also appears (in some form) in \cite{BEW08} and \cite{EMV07} 
and is a basic property of the type of automorphisms considered in this paper (also over non-$\Q$-algebras).

\vspace{0.1cm}

\begin{lemma} \label{sumcomp} 
Let $\ideaala\subset R$ be an ideal such that $\ideaala^2=(0)$. 
Suppose $G_1,G_2,H_1,H_2 \in \ideaala[X,Y]$ are given, and define $\varphi,\psi \in R[X,Y]^2$ by 
$\varphi=(X+G_1,Y+G_2)$ and $\psi=(X+H_1,Y+H_2)$. Then $\varphi\psi=\psi\varphi=(X+G_1+H_1,Y+G_2+H_2)$.

In particular, $\varphi \in \GA_2(R)$ with $\varphi^{-1}=(X-G_1,Y-G_2)$.
\end{lemma}

\begin{proof}
Straightforward.
\end{proof}

\noindent The type of tame automorphisms considered in the following proposition 
provides a foundation on which we can build many other tame automorphisms.

\vspace{0.1cm}

\begin{propo} \label{theprop}
Let $\ideaala \subseteq R$ be an ideal such that $\ideaala^2=(0)$.  
Suppose $\varphi \in \SA_2(R)$ has the form $\varphi=(X+g,Y+h)$, where $g,h \in \ideaala[X,Y]$. 
Then $\varphi\in\EA_2(R)$.
\end{propo}

\begin{proof}
Since $\ideaala^2=(0)$, $|J(\varphi)|=1+\frac{\partial g}{\partial X}+\frac{\partial h}{\partial Y}$. 
Then $\frac{\partial g}{\partial X}+\frac{\partial h}{\partial Y}=0$, and since $R$ is a $\Q$-algebra, 
this implies that there exists a polynomial $p \in \ideaala[X,Y]$ such that $g=\frac{\partial p}{\partial Y}$ 
and $h=-\frac{\partial p}{\partial X}$. Using Lemma~\ref{sumcomp}, we may assume that $p=rX^nY^m$ 
for some $r \in \ideaala,\ n,m\ge0$ and $n+m\geq1$. With Lemma~\ref{linearcomb}, 
we can write $X^nY^m$ as a $\Q$-linear combination of polynomials of the form $(X+aY)^{n+m}$, 
with $a \in \Q$. Applying Lemma~\ref{sumcomp} again, we may assume that
$$
\varphi=\left(X+kabr(X+aY)^{k-1},Y-kbr(X+aY)^{k-1}\right)
$$
where $k=n+m$ and $a,b \in \Q$. But then $\varphi=\alpha^{-1}\beta\alpha$, 
where $\alpha=(X+aY,Y)$ and $\beta=(X,Y-kbrX^{k-1})$.  Therefore $\varphi\in\EA_2(R)$.
\end{proof}

\noindent The following theorem is a special case of \cite[Theorem~4.1]{BEW08}.

\begin{theo}\label{modnil} Let $\ideaala$ be an ideal contained in 
the nilradical of $R$, and $\overline{R}=R/\ideaala$. Let $\varphi \in \SA_2(R)$. 
If $\overline{\varphi} \in \EA_2(\overline{R})$, then $\varphi \in \EA_2(R)$.
\end{theo}

\begin{proof}
Since the assumption that $\overline{\varphi} \in \EA_2(\overline{R})$ can be expressed 
using only finitely many coefficients in the ideal $\ideaala$, we may assume that $\ideaala$ is finitely generated. 
Hence it is a nilpotent ideal, say $\ideaala^m=(0)$ for some $m\geq1$. We will prove by
induction on $m$ that $\varphi$ is a composition of elementary automorphisms.

The case $m=1$ is trivial. Now suppose $m\geq2$ and let $\tilde R=R/\ideaala^{m-1}$ and $\tilde \ideaala=\ideaala/\ideaala^{m-1}$. 
Since $\tilde{\varphi}\in\SA_2(\tilde R)$, the induction hypothesis 
(applied to the ring $\tilde R$ and its ideal $\tilde \ideaala$) says that $\tilde\varphi\in\EA_2(\tilde R)$. 
Since $R\to\tilde R$ is surjective, we can lift $\tilde \varphi$ to a $\varphi_0\in\EA_2(R)$. 
Then $\varphi_0^{-1}\varphi=(X+H_1,Y+H_2)$, where $H_1,H_2 \in \ideaala^{m-1}[X,Y]$. 
The conclusion $\varphi\in\EA_2(R)$ now follows from Proposition \ref{theprop}.
\end{proof}

\begin{theo} \label{dimzerotame}
If $R$ is Artinian, then $\SA_2(R)=\EA_2(R)$.
\end{theo}

\begin{proof}
The special case of a field follows from Corollary~\ref{TAandSA} and Theorem~\ref{JvdK}. 
For the general case, let $\eta$ be the nilradical of $R$. Since $R$ is Artinian, it is well-known that 
$R/\eta$ is a product of fields. The statement now follows from Theorem~\ref{modnil} and the fact that, 
for any direct product of rings $R = R_1\times R_2$, the group $\EA_2(R)$ is canonically isomorphic 
to the direct product of groups $\EA_2(R_1)\times\EA_2(R_2)$. (And the same for $\SA_2(-)$.)
\end{proof}

\section{The square-zero principal ideal setting}
\label{specialcase}

\noindent To find the structure of the general polynomial automorphism group 
over an Artinian ring $R$, we can restrict ourselves to the case of local 
Artinian rings. 
Namely, it is well-known that $R\cong R_1\times R_2\times\cdots\times R_m$, 
a direct product of local Artinian rings. And then $\GA_2(R)$ is canonically isomorphic 
to the direct product of groups $\GA_2(R_1)\times\GA_2(R_2)\times\cdots\times\GA_2(R_m)$. 
One can readily check that this also holds if $\GA_2(-)$ is replaced by one of its mentioned subgroups.

The remainder of this paper will be focused on the case of a specific type of 
local Artinian rings, namely the ones for which the maximal ideal is principal 
and has its square equal to zero. The question of tameness over any Artinian ring 
can be reduced to this setting. We will see that the automorphism group has a 
clear structure in this case. To describe the basic aspects of this structure, 
we can use a more general setting: we suppose (for the moment) that $R$ is any ring 
containing an element $t$ satisfying $t^2=0$. In specific examples, such a ring 
is usually obtained as a factor ring of a univariate polynomial ring: $R=A[T]/(T^2)$, 
and $t=T+(T^2)$. We often use the notation ${A[t]}_2$ to denote this ring.
For this kind of ring we will give an explicit description of the group $\EA_2(R)$ 
in terms of the group $\EA_2(A)$. This will be very useful 
in the next section, when we apply this to the situation that $R$ is local Artinian.

The conjugation formulas below are crucial properties of the structure of the 
automorphism group $\SA_2(R)$. 

\begin{propo} \label{conjug}
For any $h \in R[X,Y]$ and $\alpha=(f(X,Y),g(X,Y)) \in \SA_2(R)$,
\begin{equation*} 
\alpha^{-1}(X+t\genfrac{}{}{}{1}{\partial h}{\partial Y}, Y-t\genfrac{}{}{}{1}{\partial h}{\partial X})\alpha=
(X+t\genfrac{}{}{}{1}{\partial}{\partial Y}(h(f,g)), Y-t\genfrac{}{}{}{1}{\partial}{\partial X}(h(f,g)))
\end{equation*}
In particular, if $m \in \N$ satisfies $\,m\!+\!1 \in R^*$, and
$F:=\frac{1}{m+1}f^{m+1}$ and $G:=\frac{1}{m+1}g^{m+1}$, then
$$
\alpha^{-1}(X,Y-tX^m)\alpha=(X+t\genfrac{}{}{}{1}{\partial F}{\partial
Y},Y-t\genfrac{}{}{}{1}{\partial F}{\partial X})
$$
$$
\alpha^{-1}(X+tY^m,Y)\alpha=(X+t\genfrac{}{}{}{1}{\partial G}{\partial
Y},Y-t\genfrac{}{}{}{1}{\partial G}{\partial X})
$$
\end{propo}

\begin{proof}
Let $\alpha^{-1}=(p(X,Y),q(X,Y))$. Since $t^2=0$, for any $u \in
R[X,Y]$ we have
$$
u(X+t\genfrac{}{}{}{1}{\partial h}{\partial Y}, Y-t\genfrac{}{}{}{1}{\partial h}{\partial X})=
u(X,Y)+t\genfrac{}{}{}{1}{\partial h}{\partial Y}\genfrac{}{}{}{1}{\partial u}{\partial X}-
t\genfrac{}{}{}{1}{\partial h}{\partial X}\genfrac{}{}{}{1}{\partial u}{\partial Y}=
u(X,Y)+t|\J(u,h)|
$$
Moreover, since $|\J(f,g)|=1$, the chain rule gives
$$
|\J(u(f,g),h(f,g))|=|(\J(u,h))(f,g)||\J(f,g)|=|(\J(u,h))(f,g)|
$$
The composition $\alpha^{-1}(X+t\genfrac{}{}{}{1}{\partial h}{\partial Y}, Y-t\genfrac{}{}{}{1}{\partial h}{\partial X})\alpha$ 
can now be written as
\begin{eqnarray*}
\alpha^{-1}(X+t\genfrac{}{}{}{1}{\partial h}{\partial Y}, Y-t\genfrac{}{}{}{1}{\partial h}{\partial X})\alpha & = & 
(p(X,Y)+t|\J(p,h)|, q(X,Y)+t|\J(q,h)|)
\circ(f,g)\\
 & = & (X+t|(\J(p,h))(f,g)|, Y+t|(\J(q,h))(f,g)|)\\
 & = & (X+t|\J(p(f,g),h(f,g))|, Y+t|\J(q(f,g),h(f,g))|)\\
 & = & (X+t|\J(X,h(f,g))|, Y+t|\J(Y,h(f,g))|)\\
 & = & (X+t\genfrac{}{}{}{1}{\partial}{\partial Y}(h(f,g)), Y-t\genfrac{}{}{}{1}{\partial}{\partial X}(h(f,g)))
\end{eqnarray*}
\end{proof}

\noindent These conjugation formulas naturally inspire us to make the following definition.

\vspace{0.1cm}

\begin{defi}
For any $h \in R[X,Y]$ we define $\varphi^{(h)} \in \SA_2(R)$ by 
$\varphi^{(h)}:=(X+t\genfrac{}{}{}{1}{\partial h}{\partial Y}, Y-t\genfrac{}{}{}{1}{\partial h}{\partial X})$.
\end{defi}

\vspace{0.1cm}

\begin{rema} \label{properties}
The automorphisms of the form $\varphi^{(h)}$ have the following properties:
\begin{enumerate}
\item $\varphi^{(h_1)}\varphi^{(h_2)}=\varphi^{(h_1+h_2)}$ for any $h_1,h_2 \in R[X,Y]$ 
(by Lemma~\ref{sumcomp})
\item $\alpha^{-1}\varphi^{(h)}\alpha=\varphi^{(h(f,g))}$ for $\alpha=(f,g) \in \SA_2(R)$ 
(by Proposition~\ref{conjug})
\end{enumerate}
In particular, if $m \in \N^*$ satisfies $m \in R^*$, and if $a \in R$, and $f \in R[X,Y]$ 
is one of the coordinates of an automorphism $\alpha \in \EA_2(R)$, 
then $\varphi^{(\genfrac{}{}{}{2}{a}{m}f^m)} \in \EA_2(R)$. 
Combining this with property~1. yields many tame automorphisms:
if we let $H = \genfrac{}{}{}{1}{a_1}{m_1}f_1^{m_1}+\cdots+\genfrac{}{}{}{1}{a_r}{m_r}f_1^{m_1}$, 
where $a_i \in R,\,m_i \in \N^*\cap R^*$ and each $f_i$ is a coordinate of an automorphism in $\EA_2(R)$,
then $\varphi^{(H)} \in \EA_2(R)$. In case $R={A[t]}_2$, where the ring $A$ is contained in a 
$\Q$-algebra, we have a reverse statement, displayed in Theorem~\ref{structure1}.
\end{rema}

\vspace{0.1cm}

\noindent In the proof of Theorem~\ref{structure1}, we use the following group-theoretic lemma.

\begin{lemma} \label{closure}
Let $G=H\ltimes N$ be a semidirect product of a subgroup $H$ and a normal subgroup $N$. 
Suppose we have a subset $S\subseteq N$ such that $H$ and $S$ generate the whole group $G$. 
Then $N=\,\langle h^{-1}\!sh\,:\,h \in H,\,s \in S\rangle$.
\end{lemma}

\begin{proof}
First, note that we may replace $S$ by $S \cup S^{-1}$. Now suppose $n \in N$. 
Then also $n \in G$, so we may write $n=h_1s_1\cdots h_rs_rh_{r+1}$ with $h_1,\ldots\!,h_{r+1} \in H$ 
and $s_1,\ldots\!,s_r \in S$ (some of the $h_i$ can be chosen to equal the identity). 
Viewing this $\!\!\!\!\mod\!N$, we obtain 
$1=\overline{n}=\overline{h_1}\overline{s_1}\cdots\overline{h_r}\overline{s_r}\overline{h_{r+1}}=
\overline{h_1}\cdots\overline{h_{r+1}}$, as $S\subseteq N$. 
The fact that the composition $H \hookrightarrow G \twoheadrightarrow G/N$ is an isomorphism, 
gives $1=h_1\cdots h_{r+1}$. Using this fact, we can rewrite $n$ as
$$
n = (h_1s_1h_1^{-1})((h_1h_2)s_2(h_1h_2)^{-1})\cdots((h_1\cdots h_r)s_r(h_1\cdots h_r)^{-1})
$$
\end{proof}

\noindent Before we reveal the structure of the group $\EA_2(R)$, we fix a notation for a specific 
subgroup.

\begin{defi}
$\GA_2(tR)$ denotes the subgroup of $\GA_2(R)$ consisting of those elements that have the form
$$
(X+tP(X,Y),Y+tQ(X,Y))
$$
with $P,Q \in R[X,Y]$. Furthermore, $\EA_2(tR):=\GA_2(tR)\cap\EA_2(R)$. Note that 
$\GA_2(tR)=\Ker(\GA_2(R)\to\GA_2(R/tR)) \triangleleft \GA_2(R)$. Consequently, also 
$\EA_2(tR) \triangleleft \EA_2(R)$. Obviously, if $R$ is of the form $R={A[t]}_2$, 
then $\GA_2(tR)=\GA_2(tA)$ and $\EA_2(tR)=\EA_2(tA)$. 
\end{defi}

\vspace{0.1cm}

\begin{theo} \label{structure1}
Let $A$ be a ring which is contained in a $\Q$-algebra $Q$. 
Let $R:={A[t]}_2$. 
Then, for any $\varphi_1 \in \EA_2(R)$, there exist 
a $\varphi_0 \in \EA_2(A)$ and an $H \in Q[X,Y]$ 
with $\frac{\partial H}{\partial X}, \frac{\partial H}{\partial Y} \in A[X,Y]$ such that
$$
\varphi_1 = \varphi_0\circ\varphi^{(H)} = \varphi_0 \circ (X+t\genfrac{}{}{}{1}{\partial H}{\partial
Y},Y-t\genfrac{}{}{}{1}{\partial H}{\partial X})
$$
Moreover, 
$$
H = \genfrac{}{}{}{1}{a_1}{m_1}f_1^{m_1}+\cdots+\genfrac{}{}{}{1}{a_r}{m_r}f_1^{m_1}
$$
where $a_i \in A, m_i \in \N^*$ and each $f_i$ is a coordinate of 
an automorphism in $\EA_2(A)$.
\end{theo}

\begin{proof}
Let $\varphi_1 \in \EA_2(R)$. $R=A\oplus At$, so $\EA_2(R)=\langle\EA_2(A),\EA_2(tA)\rangle$. 
Since we've also seen that $\EA_2(tA) \triangleleft \EA_2(R)$ and as it is clear that 
$\EA_2(A)\cap\EA_2(tA)=\{\textrm{id}\}$, we may conclude that $\EA_2(R)=\EA_2(A)\ltimes\EA_2(tA)$. 
So, write $\varphi_1=\varphi_0\circ\varphi_t$, with $\varphi_0 \in \EA_2(A)$ and $\varphi_t \in \EA_2(tA)$. 
Now define $S \subseteq \EA_2(tA)$ by 
$$
S = \{(X+a_itY^{m_i},Y)\,:\,a_i\in A, m_i \in \N\} \cup \{(X,Y-a_itX^{m_i})\,:\,a_i\in A, m_i \in \N\}
$$
Note that $\langle S\rangle=\{\,(X+tP(Y),Y+tQ(X))\,:\,P(Y) \in R[Y], Q(X) \in R[X]\,\}$ 
(the subgroup generated by $S$), implying that $\langle S\rangle\neq \EA_2(tA)$. For example, 
$(X+t(X-Y),Y+t(X-Y))=(X+Y,Y)(X,Y+tX)(X-Y,Y) \in \EA_2(tA)$. 
However, it is easily seen that $\EA_2(R)=\langle\EA_2(A),S\rangle$. So $G:=\EA_2(R), H:=\EA_2(A), N:=\EA_2(tA)$ and $S$ satisfy the 
requirements of Lemma~\ref{closure}. As a result, we can write
$$
\varphi_t=(\tau_1^{-1}\varepsilon_1\tau_1)(\tau_2^{-1}\varepsilon_2\tau_2)\cdots
(\tau_r^{-1}\varepsilon_r\tau_r)
$$
where each $\tau_i \in \EA_2(A)$ and each $\varepsilon_i \in S$ 
(note that $S^{-1}=S$). Then, using Proposition~\ref{conjug},
$$
\tau_i^{-1}\varepsilon_i\tau_i=
(X+a_itf_i^{m_i}\genfrac{}{}{}{1}{\partial f_i}{\partial Y},Y-a_itf_i^{m_i}\genfrac{}{}{}{1}{\partial f_i}{\partial X})=
(X+t\genfrac{}{}{}{1}{\partial h_i}{\partial Y},Y-t\genfrac{}{}{}{1}{\partial h_i}{\partial X})
$$
for some $f_i \in A[X,Y]$, and where $h_i:=\frac{a_i}{m_i+1}f_i^{m_i+1} \in Q[X,Y]$. Note that $f_i$ is a coordinate 
of an automorphism in $\EA_2(A)$. Now we can define
$H(X,Y):=h_1+\cdots+h_r$, and we derive
$$
\tau_1^{-1}\varepsilon_1\tau_1\cdots \tau_r^{-1}\varepsilon_r\tau_r=
(X+t\genfrac{}{}{}{1}{\partial H}{\partial Y},Y-t\genfrac{}{}{}{1}{\partial H}{\partial X})
$$
Obviously, $\frac{\partial H}{\partial X},\frac{\partial H}{\partial
Y} \in A[X,Y]$, whence $\varphi_1$ has the prescribed form.
\end{proof}

\noindent In case the coefficient ring is of the form ${B[t]}_2$ for a ring $B$ which is \textit{not} 
contained in a $\Q$-algebra, the above theorem can still be used to unravel the 
structure of the group $\EA_2(R)$, as is shown in Corollary~\ref{structure2}.

\vspace{0.1cm}

\begin{corol} \label{structure2}
Let $A$ be a ring which is contained in a $\Q$-algebra $Q$. 
Let $\ideaala\subseteq A$ be an ideal, and define $B:=A/\ideaala$. 
Let $R:={A[t]}_2$ and $\overline{R}:={B[t]}_2$. 
Then, for any $\varphi_1 \in \EA_2(\overline{R})$, there exist 
a $\varphi_0 \in \EA_2(B)$ and an $H \in Q[X,Y]$ 
with $\frac{\partial H}{\partial X}, \frac{\partial H}{\partial Y} \in A[X,Y]$ such that
$$
\varphi_1 = \varphi_0\circ\overline{\varphi^{(H)}} = \varphi_0 \circ (X+t\overline{\genfrac{}{}{}{1}{\partial H}{\partial
Y}},Y-t\overline{\genfrac{}{}{}{1}{\partial H}{\partial X}})
$$
Moreover, 
$$
H = \genfrac{}{}{}{1}{a_1}{m_1}f_1^{m_1}+\cdots+\genfrac{}{}{}{1}{a_r}{m_r}f_1^{m_1}
$$
where $a_i \in A, m_i \in \N^*$ and each $f_i$ is one of the coordinates of 
an automorphism in $\EA_2(A)$.
\end{corol}

\begin{proof}
Let $\varphi_1 \in \EA_2(\overline{R})$. Obviously, there exists a $\Phi_1 \in \EA_2(R)$ such 
that $\overline{\Phi_1}=\varphi_1$. The existence of $\varphi_0$ and $\varphi^{(H)}$ now follows from 
Theorem~\ref{structure1}. 
\end{proof}

\section{The case of a local Artinian ring with square-zero principal maximal ideal}

\noindent In the previous section we examined the structure of $\EA_2(R)$ in the general setting 
of a ring with a square-zero principal ideal. Now we specialize to the situation that the ring is 
local Artinian and the ideal is maximal. Whereas every automorphism over an Artinian $\Q$-algebra 
is tame (Theorem~\ref{dimzerotame}), this is not true anymore in prime characteristic, 
as is shown by the following theorem.

\begin{theo} \label{charpnontame}
Let $p$ be any prime number and $R={\F_p[t]}_2$. Then $\SA_2(R)\not\subseteq\TA_2(R)$. More
precisely, the following automorphism over $R$ is not tame:
$$
(X+tX^pY^{p-1},Y)
$$
\end{theo}

\begin{proof}
Suppose $\varphi_1:=(X+tX^pY^{p-1},Y)$ is tame. $R$ is a local ring, so 
$\varphi_1 \in \EA_2(R)$ by Lemma~\ref{primtam}. 
Now we can apply Corollary~\ref{structure2} with $A:=\Z$, $\ideaala:=p\,\Z$ and $Q:=\Q$.
Hence, there exists an $H \in \Q[X,Y]$ 
with $\frac{\partial H}{\partial X}, \frac{\partial H}{\partial Y} \in \Z[X,Y]$ such that
$$
\varphi_1 = (X+t\overline{\genfrac{}{}{}{1}{\partial H}{\partial
Y}},Y-t\overline{\genfrac{}{}{}{1}{\partial H}{\partial X}})
$$
(Note that the $\varphi_0$ of Corollary~\ref{structure2} equals the identity, since 
$\varphi_1 \in \EA_2(t\F_p)$.)
So $\overline{\frac{\partial H}{\partial Y}}=X^pY^{p-1}$, which implies
that the monomial $X^pY^p$ occurs in $H(X,Y)$, say with coefficient
$\frac{a}{b}$, where $a \in \Z\backslash\{0\}$ and $b \in \N^*$
with $\gcd(a,b)=1$. As $\frac{\partial H}{\partial Y} \in \Z[X,Y]$,
also $\frac{pa}{b}X^pY^{p-1} \in \Z[X,Y]$, whence $b\mid p$ (since $\gcd(a,b)=1$). 
Moreover, $\overline{\frac{p}{b}}\,\overline{a}=\overline{\frac{pa}{b}}=1$, so $a \notin p\,\Z$ and $b=p$.
And $\overline{\frac{\partial}{\partial X}
\frac{a}{p}X^pY^p}=\overline{a}X^{p-1}Y^p\neq0$. So the monomial $X^{p-1}Y^p$
occurs in $\overline{\frac{\partial H}{\partial X}}$, but this
contradicts the fact that $\overline{\frac{\partial H}{\partial
X}}=0$ ! (since $\varphi_1=(X+tX^pY^{p-1},Y)$) So $\varphi_1$ cannot be tame.
\end{proof}

\noindent The next example shows, that for $p=2$, a slightly modified version 
of the automorphism in Theorem~\ref{charpnontame} is tame. It is unknown to the 
author if, for all other primes $p$, the corresponding modified automorphism is tame.

\begin{exam} \label{tegenv}
Let $R={\F_2[t]}_2$. Although $(X+tX^2Y,Y) \in \SA_2(R)$ is not tame 
according to Theorem~\ref{charpnontame}, it became apparent at the end of the
proof that this is because the monomial $XY^2$ doesn't occur in the
second component of this automorphism. Then the following question
arises: is the special automorphism $(X+tX^2Y,Y-tXY^2)$ tame over $R$? Yes, it is! Since
$X^2Y$ and $XY^2$ are the partial derivatives of $\frac{1}{2}X^2Y^2$
(over $\Q$), we establish the tameness by writing this term as a linear 
combination of powers, in the style of (the end of) Remark~\ref{properties}:
\begin{eqnarray*}
\frac{1}{2}X^2Y^2 & = & \frac{1}{4}(X+Y)^4-\frac{1}{3}(Y+X^2)^3+\frac{1}{2}(Y+X^4)^2-\frac{1}{2}(Y+X^3)^2-\frac{1}{2}(X+Y^3)^2\\
 & & \!\!\!\!\!-\frac{1}{2}X^8+\frac{5}{6}X^6+\frac{1}{2}Y^6-\frac{1}{4}X^4-\frac{1}{4}Y^4+\frac{1}{3}Y^3+\frac{1}{2}X^2
\end{eqnarray*}
Applying Proposition~\ref{conjug} to each of the terms appearing in
this linear combination (taking $R={\Q[t]}_2\,$), we get
that $(X+tX^2Y,Y-tXY^2)$ equals the composition
$$
\varepsilon_0(\alpha_1^{-1}\varepsilon_1\alpha_1)(\alpha_2^{-1}\varepsilon_2\alpha_2)(\alpha_3^{-1}\varepsilon_3\alpha_3)
(\alpha_4^{-1}\varepsilon_4\alpha_4)(\alpha_5^{-1}\varepsilon_5\alpha_5)
$$
where $\varepsilon_0=(X+tY^2-tY^3+3tY^5, Y)\circ(X, Y-tX+tX^3-5tX^5+4tX^7)$ and 
$$
\begin{array}{ll}
\alpha_1=(X+Y,Y) & \varepsilon_1=(X,Y-tX^3) \\
\alpha_2=(X,Y+X^2) & \varepsilon_2=(X-tY^2,Y) \\
\alpha_3=(X,Y+X^4) & \varepsilon_3=(X+tY,Y) \\
\alpha_4=(X,Y+X^3) & \varepsilon_4=(X-tY,Y) \\
\alpha_5=(X+Y^3,Y) & \varepsilon_5=(X,Y+tX) 
\end{array}
$$
Note that this is actually a composition over ${\Z[t]}_2$. 
Viewing this composition over $R$ by calculating modulo 2, we obtain
$$
(X+tX^2Y,Y-tXY^2) \in \EA_2(R)
$$
\end{exam}

\vspace{0.5cm}

\noindent Let $p$ be a prime number. From Corollary~\ref{structure2} it follows that, 
if $R:={\Z[t]}_2$ and $\overline{R}:={\F_p[t]}_2$, then any automorphism in 
$\EA_2(\overline{R})$ is (up to an automorphism in $\EA_2(\F_p)$) of the form 
$\overline{\varphi^{(H)}}$ for some $H \in \Q[X,Y]$ with 
$\frac{\partial H}{\partial X}, \frac{\partial H}{\partial Y} \in \Z[X,Y]$. 
The automorphism $(X+tX^pY^{p-1},Y) \in \SA_2(\overline{R})$ is not of this form 
(which has been shown in the proof of Theorem~\ref{charpnontame}), so it cannot be tame. 
It is still unknown to the author whether tameness 
(more precisely: `being an element of $\EA_2(\overline{R})$') is guaranteed for all automorphisms 
over $\overline{R}$ of the form $\overline{\varphi^{(H)}}$. By Corollary~\ref{structure2}, 
this question is equivalent to the following (Question~\ref{Q2}). 
A more general version is Question~\ref{Q1}.

\vspace{0.1cm}

\begin{quest} \label{Q1}
Can every $H \in \Q[X,Y]$ with $\frac{\partial H}{\partial X}, \frac{\partial H}{\partial Y} \in \Z[X,Y]$ 
be written as a sum of the form $\genfrac{}{}{}{1}{a_1}{m_1}f_1^{m_1}+\cdots+\genfrac{}{}{}{1}{a_r}{m_r}f_1^{m_1}$, 
where $a_i \in \Z, m_i \in \N^*$ and each $f_i$ is one of the coordinates of an automorphism in $\EA_2(\Z)$ ? 
\end{quest}

\begin{quest} \label{Q2}
If the answer to Question~\ref{Q1} is negative, let $p$ be a fixed prime number. 
Does there exist, for every $H \in \Q[X,Y]$ with 
$\frac{\partial H}{\partial X}, \frac{\partial H}{\partial Y} \in \Z[X,Y]$, a sum 
$H'=\genfrac{}{}{}{1}{a_1}{m_1}f_1^{m_1}+\cdots+\genfrac{}{}{}{1}{a_r}{m_r}f_1^{m_1}$ 
as in Question~\ref{Q1}, such that $\overline{\frac{\partial H}{\partial X}}=\overline{\frac{\partial H'}{\partial X}}$ 
and $\overline{\frac{\partial H}{\partial Y}}=\overline{\frac{\partial H'}{\partial Y}}$ in $\F_p[X,Y]$ ?
\end{quest}

\vspace{0.1cm}

\noindent If Question~\ref{Q2} also has a negative answer, then the next 
challenge is to find an algorithm to decide for a given $p$ and $H$ whether 
such an $H'$ exists. Such an algorithm would thus also be an algorithm for 
tameness in $\SA_2({\F_p[t]}_2)$.

We conclude with an example of a monomial $H \in \Q[X,Y]$ which has 
the property that $\overline{\varphi^{(H)}} \in \EA_2({\F_p[t]}_2)$ for all 
primes $p\neq2$. It is unknown to the author whether this also holds for $p=2$.

\begin{exam}
It is readily verified that $\genfrac{}{}{}{1}{2}{3}X^3Y^3=\sum_{i=1}^{35} h_i$, where 
\begin{flushleft}
\begin{tabular}{@{\!}l@{\,}l@{\,}l@{\,}l@{\,}l}
 $h_1\!=\!-\genfrac{}{}{}{1}{1}{6}(X+Y)^6$ & $h_2\!=\!(Y+X^3)^4$ & $h_3\!=\!\genfrac{}{}{}{1}{5}{4}(X+Y^2)^4$ & 
 $h_4\!=\!\genfrac{}{}{}{1}{5}{4}(Y+X^2)^4$ & $h_5\!=\!-2(Y+X^6)^3$\\ [1ex]
 $h_6\!=\!-\genfrac{}{}{}{1}{5}{3}(X+Y^4)^3$ & $h_7\!=\!-\genfrac{}{}{}{1}{5}{3}(Y+X^4)^3$ & $h_8\!=\!-\genfrac{}{}{}{1}{5}{3}(X+Y^3)^3$ & 
 $h_9\!=\!-\genfrac{}{}{}{1}{5}{3}(Y+X^3)^3$ & $h_{10}\!=\!3(Y+X^{12})^2$\\ [1ex]
 $h_{11}\!=\!-2(Y+X^9)^2$ & $h_{12}\!=\!\genfrac{}{}{}{1}{5}{2}(X+Y^8)^2$ & $h_{13}\!=\!\genfrac{}{}{}{1}{5}{2}(Y+X^8)^2$ & 
 $h_{14}\!=\!\genfrac{}{}{}{1}{1}{2}(X+Y^5)^2$ & $h_{15}\!=\!\genfrac{}{}{}{1}{1}{2}(Y+X^5)^2$\\ [1ex]
 $h_{16}\!=\!-3X^{24}$ & $h_{17}\!=\!4X^{18}$ & $h_{18}\!=\!-\genfrac{}{}{}{1}{5}{2}X^{16}$ & $h_{19}\!=\!-\genfrac{}{}{}{1}{5}{2}Y^{16}$ & 
 $h_{20}\!=\!\genfrac{}{}{}{1}{2}{3}X^{12}$ \\ [1ex]
 $h_{21}\!=\!\genfrac{}{}{}{1}{5}{3}Y^{12}$ & $h_{22}\!=\!-\genfrac{}{}{}{1}{1}{2}X^{10}$ & $h_{23}\!=\!-\genfrac{}{}{}{1}{1}{2}Y^{10}$ & 
 $h_{24}\!=\!\genfrac{}{}{}{1}{5}{3}X^9$ & $h_{25}\!=\!\genfrac{}{}{}{1}{5}{3}Y^9$\\ [1ex]
 $h_{26}\!=\!-\genfrac{}{}{}{1}{5}{4}X^8$ & $h_{27}\!=\!-\genfrac{}{}{}{1}{5}{4}Y^8$ & $h_{28}\!=\!\genfrac{}{}{}{1}{1}{6}X^6$ & 
 $h_{29}\!=\!\genfrac{}{}{}{1}{1}{6}Y^6$ & $h_{30}\!=\!-\genfrac{}{}{}{1}{5}{4}X^4$\\ [1ex]
 $h_{31}\!=\!-\genfrac{}{}{}{1}{9}{4}Y^4$ & $h_{32}\!=\!\genfrac{}{}{}{1}{10}{3}X^3$ & $h_{33}\!=\!\genfrac{}{}{}{1}{16}{3}Y^3$ & 
 $h_{34}\!=\!-3X^2$ & $h_{35}\!=\!-4Y^2$
 \end{tabular}
 \end{flushleft}

\vspace{0.5cm}
 
\noindent Note that every $h_i \in \Q[X,Y]$ with 
$\frac{\partial h_i}{\partial X}, \frac{\partial h_i}{\partial Y} \in \Z[X,Y]$. 
Now, for every prime $p\neq2$ we have $\genfrac{}{}{}{1}{p+1}{3}X^3Y^3=\sum_{i=1}^{35} \genfrac{}{}{}{1}{p+1}{2}h_i$, 
from which it follows (using the same method as in Example~\ref{tegenv}) that 
$\overline{\varphi^{(\genfrac{}{}{}{3}{1}{3}X^3Y^3)}}=\overline{\varphi^{(\genfrac{}{}{}{3}{p+1}{3}X^3Y^3)}} \in \EA_2({\F_p[t]}_2)$.

\end{exam}

\vspace{1cm}

\doublespacing

\begin{center}

\end{center}

\end{document}